\documentclass[smallextended,envcountsect]{svjour3} 

\usepackage{graphicx}
\usepackage{enumerate}
\usepackage{amsmath}
\usepackage{amssymb}
\usepackage{latexsym}
\usepackage{mathrsfs}
\usepackage{tikz}
\usepackage{caption}
\usepackage{algorithm}
\usepackage{algorithmic}
\usepackage{float}
\usepackage{lipsum}
\usepackage[doipre={}]{uri} 
\newcommand{\cR}{\mathbb{R}} 
\newcommand{\cL}{\mathscr{L}} 
\newcommand{\cO}{\mathcal{O}}
\newcommand{\cF}{\mathcal{F}}
\newcommand{\cI}{\mathcal{I}}
\newcommand{\cJ}{\mathcal{J}}
\newcommand{\cK}{\mathcal{K}}
\newtheorem{assumption}{Assumption}[section]
\newtheorem{exam}{Example}[section]
\newtheorem{rem}{Remark}[section]
\newtheorem{alg}{Algorithm}[section]


\begin{document}

\title{A New Sequential Optimality Condition of Cardinality-Constrained Optimization Problems and Application}

%\subtitle{Using  the  LaTex Template}

\author{ Liping Pang  \and  Menglong Xue \and Na Xu }

\institute{Liping Pang \at
             Dalian University of Technology \\
              Dalian, China\\
              lppang@dlut.edu.cn
           \and
           Menglong Xue, Corresponding author  \at
              Dalian University of Technology \\
              Dalian, China\\
             mlx0819@mail.dlut.edu.cn
           \and 
          Na Xu \at
          Liaoning Normal University\\
          Dalian, China\\
          xuna19890223@163.com
}

\date{Received: date / Accepted: date}
%The correct dates will be entered by the editor.

\maketitle

  \begin{abstract}
  In this paper, we consider the cardinality-constrained optimization problems and propose a new sequential optimality condition for the continuous relaxation reformulation which is popular recently. It is stronger than the existing results and is still a first-order necessity condition for the cardinality constraint problems without any additional assumptions. Meanwhile, we provide a problem-tailored weaker constraint qualification, which can guarantee that new sequential conditions are Mordukhovich-type stationary points. On the other hand, we improve the theoretical results of the augmented Lagrangian algorithm. Under the same condition as the existing results, we prove that any feasible accumulation point of the iterative sequence generated by the algorithm satisfies the new sequence optimality condition. Furthermore, the algorithm can converge to the Mordukhovich-type (essentially strong) stationary point if the problem-tailored constraint qualification is satisfied.
  \end{abstract}
%\subclass{49J53 \and  49K99 \and more}

%All acknowledgements should be placed in the back of the paper after Conclusions..

   \section{Introduction}
   In recent years, $cardinality$-$constrained~optimization~problems$ (CCOP) have attracted great attention, due to its wide application in $portfolio$ \cite{protfolio_1,protfolio_2,protfolio_3}, $compressed~sensing$ \cite{compress_1}, $statistical~regression$ \cite{statistic_1,statistic_2} and other fields,  and a large number of scholars have tried to solve these problems from different perspectives. According to whether a model transformation is carried out, the existing methods are mainly divided into direct methods and indirect methods. While CCOP is non-convex and non-continuous, solving directly is extremely difficult. Therefore, this paper mainly focuses on the indirect methods, in which \cite{Schwartz_2016_siam.J} presents a new relaxed reformulation with orthogonal constraints by introducing an auxiliary variable $y$.
   
   The paper \cite{Schwartz_2016_siam.J} studied the relationship between problems the relaxation problem and CCOP, and proved that the two are equivalent in terms of global solution and feasibility. Compared with CCOP, the relaxation problem has a better structure property, such as continuity and smoothness, which allows us to have more tools to deal with the problem. But the relaxation problem is an optimization problem with orthogonal constraints, which means it is highly non-convex and difficult to solve. Because of the similarity between the relaxation problem and $mathematical~programs~with~complementarity~constraints$ (MPCC), a natural idea is to directly use MPCC's rich theoretical and numerical methods to solve the problem. However, this idea is often not feasible. For example, most of the MPCC's constraint qualification (CQ) cannot be directly applied to the relaxation problem (such as MPCC-LICQ). Even if it can be applied, it will often lead to better conclusions than MPCC. Literature \cite{Schwartz_2016_MP} remark 5.7 detailed summary of the difference between the two. This means that we cannot simply treat the relaxation problem as a special case of MPCC, but should develop problem-tailored theories and numerical algorithms. In recent years, as a large number of scholars continue to pay attention to this model, some results have been achieved.

   With the help of the $tightened~nonlinear~program$ of CCOP, denoted by $TNLP(x^{*})$. \v{C}ervinka et al. developed the classic constraint qualification to the CCOP in \cite{Schwartz_2016_MP}, proposed some CCOP customized constraint qualification (CC-CQ), and discussed the relationship between them. In addition, Kanzow et al. \cite{Schwartz_2021_ALA} adapted the quasi-normality CQ in \cite{Bertsekas_2002} and obtained a form corresponding to CCOP. And \cite{Schwartz_2021_Sequ} proposed a cone-continuity constraint qualification. At the same time, \cite{Schwartz_2016_siam.J} defines the first-order stationarity concept of the relaxation problem, called CC-Strong-stationary (CC-S-stationary)  and CC-Mordukhovich-stationary (CC-M-stationary), where CC-S-stationary is equivalent to the $Karush$-$Kuhn$-$Tucker$ (KKT) condition of the relaxation problem, and the CC-M-stationary is equivalent to the KKT condition of $TNLP(x^{*})$; \cite{Ribeiro_2020} provides a Weak-type stationarity. It is worth mentioning that, unlike CC-S-stationary and Weak-type stationarity, CC-M-stationary is only related to the original variable $x$, and \cite{Schwartz_2021_ALA} proves that CC-S-stationary and CC-M-stationary are equivalence in the original variable space. Consequently, this paper will focus on CC-M-stationarity.
	 	 
   Because of the similarity between the relaxation problem and MPCC, some researchers try to apply the classic algorithm of MPCC to solve the relaxation problem. \cite{Schwartz_2016_siam.J} and \cite{Schwartz_2018} respectively applied two classic MPCC's regularization methods to the relaxation problem, and both obtained better convergence than general MPCC. However, the regularization strategy is actually to further relax the relaxation problem into a sequence of regular subproblems and obtain the solution of the relaxation problem by solving the regular subproblems. Can the relaxation problem be solved directly without further relaxation? This is an issue worthy of attention. In this paper, we have made a great answer to this. The main contributions of this paper are as follows:
	 
   \begin{enumerate}[$\bullet$]
	  \item {We propose a new sequential optimality condition: CC-PAM-stationarity. In recent years, the application of sequential optimality condition to stop criteria and uniform convergence analysis of algorithms has received great attention. In this area, several sequential optimality conditions have been proposed for $nonlinear~programming$  (NLP) \cite{Andreani_2011_AKKT,Andreani_2010_CAKKT,Haeser_2011_SAKKT,Andreani_2019_PAKKT}, where \cite{Andreani_2019_PAKKT} gives the relationship between them. However, there are still very few relevant results about CCOP. \cite{Ribeiro_2020} establishes a sequential optimality condition, called CC-approximate weak stationarity (CC-AW-stationarity), but this condition is based on the $(x,y)$ space. Therefore, Kanzow et al. \cite{Schwartz_2021_Sequ} proposed CC-approximate Mordukhovich stationarity (CC-AM-stationarity), which is only related to $x$, and a proof is given that it is equivalent to CC-AW-stationary. However, for some problems, the number of CC-AM-stationary points is numerous, and these points are often far from the optimal solutions (e.g. Example \ref{exam1}). In order to obtain fewer optimal candidate points, we propose CC-PAM-stationarity, which is strictly stronger than CC-AM-stationarity, and we show that it is a necessary condition of CCOP without any assumptions.}
	  
	  \item {We define a new problem-tailored constraint qualification, called CC-PAM-regularity, which is weaker than CC-AM-regularity proposed in \cite{Schwartz_2021_Sequ}. We prove that any CC-M-stationary point is CC-PAM-stationary, and conversely, CC-PAM-stationary point is CC-M-stationary if CC-PAM-regularity condition is satisfied. In other words, CC-PAM-regularity condition is a CC-CQ. Borrowing the notation of reference \cite{ALGENCAN}, this constraint qualification is called strict constraint qualification (SCQ). Furthermore, we show that CC-PAM-regularity condition is the weakest SCQ relative to CC-PAM-stationarity.}
	  
	  \item{We apply CC-PAM-stationarity to safeguarded augmented Lagrangian method and further improve its convergence. Different from the regularization methods, the literature \cite{Schwartz_2021_ALA} and \cite{Schwartz_2021_Sequ} try to directly apply the safeguarded augmented Lagrangian method of the general NLP to the relaxation problem and \cite{Schwartz_2021_ALA} uses the corresponding solver ALGENCAN\cite{ALGENCAN,ALGENCAN_1,ALGENCAN_2} to solve the $portfolio$ problem verify the advantages of the augmented Lagrangian algorithm over the regularization methods. In addition, the above two regularization methods both require accurate KKT points for their subproblems, while the safeguarded augmented Lagrangian method only requires subproblems to be solved inaccurately. Kanzow et al. \cite{Schwartz_2021_Sequ} show that any feasible limit point of safeguarded augmented Lagrangian method is CC-AM-stationary. And we proved that under mild conditions such as semialgebraic properties (or the same conditions as \cite{Schwartz_2021_Sequ}), these points are CC-PAM-stationary, which is strictly better than CC-AM-stationary. If additional conditions of CC-PAM-regularity hold, they will be CC-M-stationary points.}
	
  \end{enumerate}
 
	 The organization is as follows: we give some basic definitions and preliminary conclusions in Sect.2; propose a new sequential optimality condition in Sect.3, and defines a new problem-tailored constraint qualification in Sect.4. The convergence of safeguarded augmented Lagrangian method is discussed in Sect.5, and Sect.6 is a simple summary.
	 
	 Notation: $I_{g}(x)=\{i:g_{i}(x)=0,~i=1,\dots,m\}$, $I_0(x)=\{\imath:x_\imath=0,~\imath=1,\dots,n\}$, $I_\pm(x)=\{\imath:x_\imath\neq 0,~\imath=1,\dots,n\}$, $x_+=\max\{x,0\}$, $|\cdot|$ denotes $l_1$-$norm$, $\parallel\cdot\parallel$ is Euclidean norm, $\parallel\cdot\parallel_\infty$ denote infinity norm, $\parallel \cdot\parallel_0$ denotes $l_0$ norm (the number of non-zero elements), $e=(1,\dots,1)^T\in\cR^{n}$, $e_i$ is a vector where only the ith component is 1 and all the others are 0, $x\circ y$ denotes Hadamard product of $x$ and $y$.
	    
\section{Preliminaries}
   In this paper, we consider the optimization problems
	\begin{equation}
      \min~ f(x) \quad s.t.\hspace{1em}g(x)\leq 0,
      \quad h(x)=0,\quad \parallel x\parallel_{0}\ \leq \kappa,
      \label{primary_problem}
    \end{equation}
   where $\kappa$ is an integer and $\kappa<n$, $f: \cR^{n}\rightarrow\cR$, $g: \cR^{n}\rightarrow\cR^{m}$, $h: \cR^{n}\rightarrow\cR^{p}$ is continuously differentiable, and $\parallel x\parallel_0$ is also called cardinality of $x$. Thus, the problem \eqref{primary_problem} is called a $cardinality$-$constrained~optimization~problems$ (CCOP). Let $x^*\in\cR^n$, the tightened NLP problem ($TNLP(x^*)$) of CCOP defined as
   	   \begin{equation}
   	   	 \min~ f(x) \quad s.t.\hspace{1em}g(x)\leq 0,
   	   	 \quad h(x)=0,\quad x_\imath=0,\hspace{0.5em} \imath\in I_{0}(x^{*}).
   	   	 \label{TNLP}
   	   \end{equation}   
   And the relaxation problem of CCOP is defined as 
   \begin{equation}
         \left\{\begin{aligned}
          &\min\quad f(x)\\
           &~s.t.\hspace{1em}g(x)\leq 0,\quad h(x)=0,\\
           &\qquad x\circ y=0,\quad n-\kappa-e^{T}y\leq 0,\quad y\leq e.
          \end{aligned}\right.\label{relax_problem}
   \end{equation}
    Note that the problem \eqref{relax_problem} is one less non-negative constraint than the form in \cite{Schwartz_2016_siam.J}. The literature \cite{Schwartz_2021_ALA} shows that this change will not affect the original conclusion and can lead to a larger feasible set. In the introduction, we have mentioned the relationship between CCOP and the problem \eqref{relax_problem}. Below we will give specific conclusions.  
    \begin{proposition}
    \textup{\cite{Schwartz_2016_siam.J}}~Let $x\in\cR^n$, then the following statements hold.
    		\begin{itemize}
    			\item {If $x$ is a feasible point (or golbal minimizer) of CCOP if and only if there exists $y\in\cR^n$ such that $(x,y)$ is feasible point (or golbal minimizer) of the problem \eqref{relax_problem}.}
    			\item{If $x$ is a local minimizer of CCOP, then there exists $y\in\cR^n$ such that $(x,y)$ is local minimizer of the problem \eqref{relax_problem}; Conversely, if $(x,y)$ is a local minimizer of the problem \eqref{relax_problem} and $\|x\|_0=\kappa$ holds, then $x$ is a local minimizer of CCOP.}
    		\end{itemize}
    \end{proposition}
    
    There are several stationarity concepts with the relaxation problem \eqref{relax_problem}. 
	\begin{definition}
		\cite{Schwartz_2016_siam.J}~Let $(x^{*},y^*)$ be feasible for \eqref{relax_problem}, then it is called  
		 \begin{enumerate}[(1)]
		 \item {CC-S-stationary, if there exists $\{(\lambda,\mu,\gamma)\}\in \cR^m\times\cR^p\times\cR^n$ such that 
		 \begin{enumerate}[$\bullet$]
		 	\item {$\nabla f(x^{*})+\nabla g(x^{*})\lambda+\nabla h(x^{*})\mu+\gamma= 0$;}
		 	\item {$\lambda_{i}=0,~\forall i\notin I_{g}(x^{*});~\gamma_{\imath}=0,~for~all~ \imath ~such~that ~y_\imath^*=0$.}
		 \end{enumerate}}		
		 \item{CC-M-stationary, if there exists $\{(\lambda,\mu,\gamma)\}\in \cR^m\times\cR^p\times\cR^n$ such that
		       \begin{enumerate}[$\bullet$]
		     		\item {$\nabla f(x^{*})+\nabla g(x^{*})\lambda+\nabla h(x^{*})\mu+\gamma= 0$;}
		     		\item {$\lambda_{i}=0,~\forall i\notin I_{g}(x^{*});~\gamma_{\imath}=0,~\forall~ \imath\in I_{\pm}(x^{*})$.}
		     	\end{enumerate}}  	
		 \end{enumerate}				
	\end{definition}
	
	Obviously, CC-M-stationarity is weaker than CC-S-stationarity, but it only depends on the variable $x$, which can be used as the optimality measure of CCOP. Another important reason to focus on CC-M-stationary points in this paper is because of the validity of the following conclusion.
    \begin{proposition}
    	\label{prop2}
    	\textup{\cite{Schwartz_2021_ALA}}~Let $(x,y)$ is feasible for \eqref{relax_problem}, if $(x,y)$ is a CC-M-stationary point, then there exists $z\in\cR^n$ such that $(x,z)$ is a CC-S-stationary point.
    \end{proposition}    
    
    Let us now recall a basic concepts that needed for theoretical analysis \cite{Variational}. The upper limit of set-valued maps $\Theta:~\cR^n\rightrightarrows\cR^m$ is 
	\[\limsup_{x\rightarrow x^*}\Theta(x):=\left\{z:\exists x^k\rightarrow x^*,~\exists z^k\rightarrow z,~z^k\in \Theta(x^k)\right\}.\]
   	
	For a function $l:\cR^n\rightarrow\cR$, the (lower) level set $l_\alpha$ is defined as
    \[l_\alpha:=\left\{~x\in\cR^n:~l(x)\leq \alpha\right\}.\]
    If any level set of the function $l$ is bounded, then the function $l$ is said to be level bounded. Since some of the conclusions of this article are obtained under the assumption of semialgebraic, let us briefly introduce the basic definition and properties of semialgebraic. We say the set $C\subseteq\cR^n$ is semialgebraic if it can be written as a finite
    union of sets of the form
	\[\left\{x \in \cR^{n}: u_{i}(x)=0,~v_{i}(x)<0,~i=1,\ldots, p\right\},\]
	where $u_i(x)$, $v_i(x)$ are polynomial functions. A function is called semialgebraic if its graph is a semialgebraic set, obviously polynomial functions are semialgebraic. Because of their strong stability, semialgebraic functions are a very broad class of functions.
    \begin{lemma}
		\textup{\cite{Attouch_2010}}~The following properties hold.
		\begin{itemize}
			\item {Linear combination of finite number of semialgebraic functions is semialgebraic.}
			\item {Composition of semialgebraic functions is semialgebraic.}
			\item{Generalized inverse of a semi-algebraic function is semialgebraic.}
			\item{Let $F(x)=\sup\limits_{y\in C}f(x,y)$ and $G(x)=\inf\limits_{y\in C}f(x,y)$, if the set $C$ and function $f$ are semialgebraic, then both F and G are semialgebraic.}
		\end{itemize}
	\label{semialge}
	\end{lemma}	
	
	Semialgebraic functions have another important property, they satisfy the Kurdyka-\L ojasiewicz property.
	\begin{definition}[KL property]
		\textup{\cite{Attouch_2009}}~We say the function $f$ satisfy the Kurdyka-\L ojasiewicz property, if for any limiting-critical point $x^*$ ($0\in \partial f(x^*)$), there exist $\epsilon,~C>0$, $\theta\in[0,1)$ such that 
		\begin{equation}
		\label{KL}
		C\hspace{0.1em}|f(x)-f(x^*)|^\theta\leq\|v\|,\quad \forall \hspace{0.1em}\|x-x^*\|\leq\epsilon,~v\in\partial f(x^*).
		\end{equation}
	\end{definition}    
	
   And if the constraint set of NLP is denoted as
   	\[X:=\{x:~g_{i}(x)\leq 0,~i=1,\dots,m;~h_{j}(x)=0,~ j=1,\dots ,p\},\]
   	then standard MFCQ is defined as follows.
   	\begin{definition}[MFCQ]
   	   		\label{MFCQ}
   	   		Let $x\in X$, then we say $x$ satisfies Mangasarian-Fromovitz CQ, if the gradient vectors $\nabla h_j(x)~(j=1,\dots,p)$ are linearly independent, and there exists $d\in\cR^n$ such that 
   	   				\[\nabla h_j(x)^Td=0,~\forall j=1,\dots,p,\qquad \nabla g_{i}(x)^Td<0,~\forall i\in I_g(x).\]
   	\end{definition}
  	
\section{A New Sequential Optimality Condition}  	
	Recently, due to its excellent properties, sequential optimality conditions are very popular. Although there have been many theoretical results on NLP problems, to avoid auxiliary variables, we did not directly apply the results of NLP to problem \eqref{relax_problem} but proposed new problem-tailed sequential optimality conditions.
   		\begin{definition}[CC-PAM-stationary]
   			Let $x^{*}$ is feasible for CCOP, we say that $x^{*}$ is CC-positive approximate Mordukhovich stationary, if there exist $\{(x^{k},\lambda^{k},\mu^{k},\gamma^{k})\}\in\cR^n\times\cR_+^m\times\cR^p\times\cR^n$ such that:
   			\begin{enumerate}[$(a)$]
   				\item {$x^{k}\rightarrow x^{*}$, $\nabla f(x^{k})+\nabla g(x^{k})\lambda^{k}+\nabla h(x^{k})\mu^{k}+\gamma^{k}\rightarrow 0$;}
   				\item {$\lambda_{i}^{k}=0,~\forall i\notin I_{g}(x^{*}),~\gamma_{\imath}^{k}=0,~\forall \imath\in I_{\pm}(x^{*})$;}
   				\item {$\lambda_{i}^{k}g_{i}(x^{k})>0$, if $\lim\limits_{k}\frac{\lambda_{i}^{k}}{\pi_k}>0$;}
   				\item {$\mu_{j}^{k}h_{j}(x^{k})>0$, if $\lim\limits_{k}\frac{|\mu_{j}^{k}|}{\pi_k}>0$;}
   				\item {$\gamma_{\imath}^{k}x_{\imath}^{k}>0$, if $\lim\limits_{k}\frac{|\gamma_{\imath}^{k}|}{\pi_k}>0$;}
   			\end{enumerate}
   		where $\pi_k=\parallel (1,\lambda^{k},\mu^{k},\gamma^{k})\parallel_{\infty}$, the sequence that satisfy the conditions $(a)$-$(e)$ are called a CC-PAM sequence.
   			\label{CC-PAM}	
   		\end{definition}
   		
   	The condition in Definition \ref{CC-PAM} only needs to be true for sufficiently large k. For example, if there is $N$, the conditions $(a)$-$(e)$ are satisfied when $k\geq N$, then you can set $\hat{x}^{k}=x^{N+k}$, and the new sequence obtained is the CC-PAM sequence. Observe that, the conditions $(a)$-$(b)$ are the same as CC-AM-stationarity, so there are the following conclusions.
   		\begin{proposition}
   			Let $x^{*}$is feasible for CCOP, if $x^{*}$ is a CC-PAM-stationary point, then it's a CC-AM-stationary point.
   		\end{proposition}
   		
   	The converse of the above conclusion is untenable, as shown in the following example.
   		\begin{exam}
   			\label{exam1}
   			We consider 
   			\begin{equation}
   			\label{exam1_pro}
   			\min_{x\in\cR^3}~\frac{1}{2}\left[(x_1-1)^2+(x_2-1)^2\right]\quad s.t.~x_1x_3\leq0,~\|x\|_0\leq 2.
   			\end{equation}
   			Obviously, the problem \eqref{exam1_pro} has the only global optimal solution $(1,1,0)^T$. Let $x=(a,1,0)^T$, where $0<a<1$. Take 
   			\[x^k=(a,1,\frac{1-a}{k}),\quad\lambda^k=k,\quad\gamma^k=(0,0,-ka)^T.\]
   			It is easy to verify that the above sequence satisfies the conditions $(a)$-$(b)$, that is, $x$ is a CC-AM-stationary point; but it is not CC-PAM-stationary, because the sequence meets the conditions $(a)$-$(b)$, $\gamma_3^k$ and $x_3^k$ must have different signs, which violates the condition $(e)$. 
   		\end{exam}
   		
   	As can be seen from Example \ref{exam1}, the number of CC-AM-stationary points is numerous, and these points are far from the optimal solution. While CC-PAM-stationary points contain fewer candidate points, that is, it is strictly superior to CC-AM-stationarity. The following Theorem \ref{thm1} states that CC-PAM-stationarity is a necessary optimality condition for CCOP without any additional assumptions.
   		
   		\begin{theorem}
   				\label{thm1}
   				Let $x^{*}$ is a local minimizer of CCOP, then $x^{*}$ is a CC-PAM-stationary point.
   		\end{theorem}
   		{\it Proof} 
   		 %\begin{proof}
   		 ~If $x^{*}$ is a local minimizer of CCOP, then it is also a local minimizer of $TNLP(x^{*})$, there exist $\epsilon>0$ such that $x^{*}$ is the only global minimizer for the following problem
   				\begin{equation}
   				\left\{\begin{aligned}
   				\min\quad &f(x)+\frac{1}{2}\parallel x-x^{*}\parallel^{2}&\\
   				~s.t.\hspace{1em}&g(x)\leq 0,\quad h(x)=0,\quad \\
   				&x_\imath=0,\hspace{0.5em} \imath\in I_{0}(x^{*}), \quad \parallel x-x^{*}\parallel\leq \epsilon.
   				\end{aligned}\right.\label{local_TNLP}
   				\end{equation}
   				Let $p(x)=\parallel h(x)\parallel^2+\parallel g(x)_+\parallel^2+\sum\limits_{\imath\in I_0(x^*)}x_\imath^2$, we define the local penalized problem
   				\begin{equation}
   				\left\{\begin{aligned}
   				\min\quad &f(x)+\frac{1}{2}\parallel x-x^{*}\parallel^{2}+\frac{M_k}{2}p(x)&\\
   				~s.t.\hspace{1em}&\parallel x-x^{*}\parallel\leq\epsilon,
   				\end{aligned}\right.\label{local_penality_TNLP}
   				\end{equation}
   				where $0<M_k\rightarrow+\infty$. For all $M_k$, the objective function of the problem \eqref{local_penality_TNLP} is continuous and the feasible set is compact, there must exist a global minimizer, denoted as $x^k$. Meanwhile, the sequence $\{x^k\}$ is bounded, there must be a convergent subsequence. For simplicity, let us set $x^k\rightarrow\bar{x}$. The following proves that $\bar{x}=x^*$.
   				
   				Since $x^k$ is the global minimizer of the problem \eqref{local_penality_TNLP}, then
   				\begin{equation}
   				f(x^k)+\frac{1}{2}\parallel x^k-x^{*}\parallel^{2}+\frac{M_k}{2}p(x^k)\leq f(x^*).
   				\label{penality_1}
   				\end{equation}
   				Divide both sides of \eqref{penality_1} by $M_k$ and take the limit, we obtain $p(\bar{x})\leq0$. So $\bar{x}$ is feasible for the local problem \eqref{local_TNLP}. In addition, from \eqref{penality_1}
   				\[f(x^k)+\frac{1}{2}\parallel x^k-x^*\parallel^2\leq f(x^*).\]
   				Letting $k\rightarrow+\infty$ yields
   				\[f(\bar{x})+\frac{1}{2}\parallel \bar{x}-x^*\parallel^2\leq f(x^*),\]
   				but $x^*$ is the only global minimum point of the problem \eqref{local_TNLP}, there must be $\bar{x}=x^*$, that is, $x^k\rightarrow x^*$.
   				
   				When $k$ is sufficiently large, there is obviously $\parallel x^k-x^*\parallel\leq\epsilon$. For simplicity, let's set $\{x^k\}\subseteq \{x:\parallel x^k-x^*\parallel\leq\epsilon\}$. From the necessary optimality condition of the problem \eqref{local_penality_TNLP} we obtain
   				\[\nabla f(x^k)+\nabla g(x^k)(M_kg(x^k)_+)+\nabla h(x^k)(M_kh(x^k))+\sum_{\imath\in I_0(x^*)}M_kx_{\imath}^ke_\imath=x^*-x^k.\]
   				And we define
   				\[\lambda^k=M_kg(x^k)_+,\quad \mu^k=M_kh(x^k),\quad \gamma_{\imath}^k=M_kx_{\imath}^k,\imath\in I_0(x^*),\quad \gamma_{\imath}^k=0,\imath\in I_\pm(x^*),\]
   				then
   				\[\nabla f(x^k)+\nabla g(x^k)\lambda^k+\nabla h(x^k)\mu^k+\gamma^k\rightarrow 0.\]
   				
   				 For all $i\notin I_{g}(x^*)$, there is $g_i(x^*)<0$, then $g_i(x^k)<0$ when $k$ is sufficiently large, so $\lambda_{i}^k=0$. Meanwhile, by the definition of $\gamma^k$, obviously
   				\[\gamma_{\imath}^k=0,\quad\forall\imath\in I_\pm(x^*).\]
   				
   				If $\lambda_{i}^k > 0$, since $\lambda_i^k=M_kg_i(x^k)_+$ then $g_i(x^k)>0$, so
   				\[\lambda_i^kg_i(x^k)>0.\]
                Similarly, if $\mu_{j}^k\neq 0$, then $h_j(x^k)\neq 0$; and if $\gamma_\imath^k\neq 0$, we obtain $x_\imath^k\neq0$. So
	            \begin{equation*}
   				\mu_{j}^kh_j(x^k)=M_k(h_j(x^k))^2>0~~and~~\gamma_\imath^kx_\imath^k=M_k(x_\imath^k)^2>0.
   				\end{equation*}   				

   				In summary, $\{(x^k,~\lambda^k,~\mu^k,~\gamma^k)\}$ is a CC-PAM sequence, that is, $x^*$ is a CC-PAM-stationary point.
   		\qed
   		%\end{proof}
   		
   	Theorem \ref{thm1} and Example \ref{exam1} show that CC-PAM-stationary point we proposed can be used as a candidate point for the optimal solution, and it is more suitable as a measure of optimality than CC-AM-stationarity. On the other hand, another advantage of the sequential optimality condition is that it has nothing to do with the specific algorithm. That is, CC-PAM-stationary point has some theoretical properties, and any algorithm that can generate the CC-PAM-stationary point also has the same nature. Therefore, the existence of sequential optimality conditions provides a tool for establishing a unified framework for optimality theory.
	
	The converse of the above conclusion is untenable, as shown in the following example.
	 \begin{exam}
	 	In two-dimensional space, consider a simple geometric problem, set $z=(1,1)^T$, find the point closest to $z$ on the coordinate axis. This problem can be modeled as
	 	\begin{equation}
	 	\label{exam2_pro}
	 	\min~\frac{1}{2}\left[(x_1-1)^2+(x_2-1)^2\right]\quad s.t.~\|x\|_0\leq 1.
	 	\end{equation}
	 	Obviously, $(1,0)^T$ and $(0,1)^T$ are the two global optimal solutions of the problem. The following shows that $x^*=(0,0)^T$ is a CC-PAM-stationary point. Take
	 	\[x^k=(\frac{1}{k+1},\frac{1}{k+1})^T\quad \gamma^k=(1-\frac{1}{k+1},1-\frac{1}{k+1})^T,\]
	 	it is easy to verify that $\{(x^k,\gamma^k)\}$ satisfies the conditions $(a)$-$(e)$, that is, $x^*=(0,0)^T$ is a CC-PAM-stationary point, but it is not a local minimizer of the problem \eqref{exam2_pro}.
	 \end{exam}
	 
	 We know that CC-M-stationarity is stronger than the CC-AM-stationarity, and CC-PAM-stationarity we proposed is also strictly better than CC-AM-stationarity. An interesting question is whether CC-M-stationarity is stronger than CC-PAM-stationarity, and under what conditions are the two equivalent. This issue will be described in detail in the next section.  	
   	
\section{A New Constraint Qualification}

   Let $x^*$ be feasible for CCOP, $\alpha\geq 0$, $\beta\geq 0$, $x\in\cR^n$, we defined the set
	\begin{equation*}
	\Theta(x,\alpha,\beta)=\left\{\nabla g(x)\lambda+\nabla h(x)\mu+\gamma\left|
	\begin{aligned}
	(\lambda,\mu,\gamma)&\in\cR_+^m\times\cR^p\times\cR^n,\\
	\lambda_{i}=0,~\forall i\notin I_g&(x^*),~\gamma_{\imath}=0,~~\forall\imath\in I_\pm(x^*),\\
	\lambda_{i}g_i(x)\geq\alpha,~~&if~~\lambda_{i}~>\beta\parallel(1,\lambda,\mu,\gamma)\parallel_\infty,\\
	\mu_{j}h_j(x)\geq\alpha,~~&if~|\mu_{j}|>\beta\parallel(1,\lambda,\mu,\gamma)\parallel_\infty,\\
	\gamma_{\imath}x_\imath\geq\alpha,~~&if~|\gamma_{\imath}|>\beta\parallel(1,\lambda,\mu,\gamma)\parallel_\infty\\
	\end{aligned}\right.\right\}.
	\end{equation*}
	Obviously, if $x^*$ is a CC-M-stationary point, then it can be written as
	\begin{equation}
	-\nabla f(x^*)\in\Theta(x^*,0,0).
	\label{CC_M_reform}
	\end{equation}
	
	In addition, CC-PAM-stationarity can be expressed as the limit form of the set sequence, and the following conclusions are established.
	\begin{lemma}
		$x^*$ is CC-PAM-stationary $\iff -\nabla f(x^*)\in\limsup\limits_{
			x^k\rightarrow x^*,
			\alpha\downarrow 0,\beta\downarrow 0
		}\Theta(x,\alpha,\beta).$
	\label{CC-PAM-set}
	\end{lemma}
	{\it Proof}~~
	%\begin{proof}
	"$\Rightarrow$" ~If $x^*$ is CC-PAM-stationary, from Definition \ref{CC-PAM}, there is a sequence $\{(x^k,\lambda^k,\mu^k,\gamma^k)\}$ such that the conditions $(a)$-$(e)$ holds. Take
		\[\theta^k=\nabla g(x^k)\lambda^k+\nabla h(x^k)\mu^k+\gamma^k,\]
		we know $\nabla f(x^k)+\theta^k\rightarrow 0$, by the condition $(a)$, then $\theta^k\rightarrow\theta^*=-\nabla f(x^*)$. 
		
		To prove the conclusion, just find the appropriate $\{\alpha_k\}$, $\{\beta_k\}$ such that $\alpha_k\downarrow 0$, $\beta_k\downarrow 0$ and $\theta^k\in\Theta(x^k,\alpha_k,\beta_k)$.
		
		Let $\pi_k=\parallel (1,\lambda^k,\mu^k,\gamma^k)\parallel_\infty$, then the sequence $\left\{\frac{(\lambda^k,\mu^k,\gamma^k)}{\pi_k}\right\}$ must have a convergent subsequence. For simplicity, let's set it to converge. Let  
		\begin{equation*}
		I=\left\{i\left|\lim\limits_{k\rightarrow\infty}\frac{\lambda_{i}^k}{\pi_k}>0\right.\right\},~~J=\left\{j\left|\lim\limits_{k\rightarrow\infty}\frac{|\mu_{j}^k|}{\pi_k}>0\right.\right\},~~K=\left\{\imath\left|\lim\limits_{k\rightarrow\infty}\frac{|\gamma_{\imath}^k|}{\pi_k}>0\right.\right\},
		\end{equation*}
		so we can take
		\[\alpha_k=\min\left\{(\lambda_{i}^k g_{i}(x^k))_{i\in I},~(\mu_{j}^k h_{j}(x^k))_{j\in J},~(\gamma_{\imath}^k x_\imath^k)_{\imath\in K},~\frac{1}{k}\right\}.\]
		Obviously there is $\alpha_k\rightarrow 0$, and when k is sufficiently large, we obtain
		\[\frac{\lambda_{i}^k}{\pi_k}>\max\left\{(\frac{\lambda_{i}^k}{\pi_k})_{i\notin I},~(\frac{|\mu_{j}^k|}{\pi_k})_{j\notin J},~(\frac{|\gamma_{\imath}^k|}{\pi_k})_{\imath\notin K},~\frac{1}{k}\right\}\quad \forall i\in I.\]
		Regarding $\mu$, $\gamma$ has a similar conclusion. And let
		\[\beta_k=\max\left\{(\frac{\lambda_{i}^k}{\pi_k})_{i\notin I},~(\frac{|\mu_{j}^k|}{\pi_k})_{j\notin J},~(\frac{|\gamma_{\imath}^k|}{\pi_k})_{\imath\notin K},~\frac{1}{k}\right\},\]
		Obviously there is $\beta_k\rightarrow 0$. Combining the non-negativity of $\{\alpha_k\}$ and $\{\beta_k\}$, we can set $\alpha_k\downarrow 0$ and $\beta_k\downarrow 0$ (if necessary, subsequence can be taken). And we obtain
		\[\theta^k\in\Theta(x^k,\alpha_k,\beta_k).\]
		"$\Leftarrow$"~~By hypothesis, there is $\{x^k\}$, $\{\alpha_k\}$, $\{\beta_k\}$ such that
		\begin{equation*}
		x^k\rightarrow 0,~\alpha_k\downarrow 0,~\beta_k\downarrow 0,~\Theta(x^k,\alpha_k,\beta_k)\ni\theta^k\rightarrow-\nabla f(x^*).
		%\label{theta_k}
		\end{equation*}
		Therefore, there is $\{(\lambda^k,\mu^k,\gamma^k)\}$ such that 
		\[\theta^k=\nabla g(x^k)\lambda^k+\nabla h(x^k)\mu^k+\gamma^k,\]
		so $\nabla f(x^k)+\theta^k\rightarrow 0$, that is
		\[\nabla f(x^k)+\nabla g(x^k)\lambda^k+\nabla h(x^k)\mu^k+\gamma^k\rightarrow 0.\]
		
		Since $\theta^k\in\Theta(x^k,\alpha_k,\beta_k)$, we obtain $\lambda_{i}^k=0~(\forall i\notin I_g(x^*)),~\gamma_{\imath}^k=0~(\forall\imath\in I_\pm(x^*))$. In addition, if $\lim\limits_{k\rightarrow\infty}\frac{|\mu_{j}^k|}{\pi_k}>0$, then $\frac{|\mu_{j}^k|}{\pi_k}>\beta_k$ for all $k$ large enough, this implies
		\[\mu_{j}^kh_{j}(x^k)\geq\alpha_k>0,\]
		that is the condition (b) is satisfied. There are similar conclusions about $\lambda$ and $\gamma$, so $x^{*}$ is CC-PAM-stationary.		
	\qed
   	%\end{proof}
   	
     Now, we give a new regularity condition.
	\begin{definition}[CC-PAM-regularity]
		We say a feasible point $x^*$ of CCOP satisfies the CC-PAM-regularity condition, if 
		\[\limsup\limits_{
			x\rightarrow x^*,~
			\alpha\downarrow 0,~\beta\downarrow 0
		}\Theta(x,\alpha,\beta)\subseteq\Theta(x^*,0,0).\]
	\label{CC_PAM_regular}
	\end{definition}
    \begin{rem}
   CC-PAM-regularity condition is weaker than CC-AM-regularity condition. Take
    \begin{equation*}
    \hat{\Theta}(x)=\left\{\nabla g(x)\lambda+\nabla h(x)\mu+\gamma\left|
    \begin{aligned}
    (\lambda,\mu,\gamma)&\in\cR_+^m\times\cR^p\times\cR^n,\\
    \lambda_{i}=0,~\forall i\notin I_g&(x^*),~\gamma_{\imath}=0,~~\forall\imath\in I_\pm(x^*)
    \end{aligned}\right.\right\}.
    \end{equation*}	
    Obviously there is $\hat{\Theta}(x^*)=\Theta(x^*,0,0)$, and for all $\alpha\geq 0$, $\beta\geq 0$, $\forall x\in\cR^n$ we have $\Theta(x,\alpha,\beta)\subseteq\hat{\Theta}(x)$. Therefore, if CC-AM-regularity condition is established, that is, $\limsup_{x\rightarrow x^*}\hat{\Theta}(x)\subseteq\hat{\Theta}(x^*)=\Theta(x^*,0,0)$, we obtain $\limsup_{x\rightarrow x^*,~
    	\alpha\downarrow 0,~\beta\downarrow 0}\Theta(x,\alpha,\beta)\subseteq\Theta(x^*,0,0)$.
	\end{rem}
	
	Now we give the relationship between CC-PAM-stationarity and CC-M-stationarity.
	\begin{theorem}
		\label{PAM->M}
		Let $x^*$ is feasible for CCOP, the following statements hold.
		\begin{enumerate}[$(\romannumeral1)$]
			\item {If $x^*$ is a CC-M-stationary point, then it is a CC-PAM-stationary point.}
			\item {If $x^*$ is a CC-PAM-stationary point, and CC-PAM-regularity condition holds at $x^*$, then $x^*$ is a CC-M-stationary point.}
			\item {If for any continuous differentiable function $f$, the following relationship holds:
				\[x^*~\text{CC-PAM-}stationary\quad \Longrightarrow\quad x^*~\text{CC-M-}stationary,\]
				then CC-PAM-regularity condition holds at $x^*$.}
		\end{enumerate}
	\end{theorem}
	{\it Proof}~$(\romannumeral1)$
	%\begin{proof}$(\romannumeral1)$
	~If $x^*$ is a CC-M-stationary point, then 
			\[-\nabla f(x^*)\in\Theta(x^*,0,0)\subseteq\limsup\limits_{
				x\rightarrow x^*,~
				\alpha\downarrow 0,~\beta\downarrow 0
			}\Theta(x,\alpha,\beta).\]
			
		$(\romannumeral2)$~By Lemma \ref{CC-PAM-set}, we obtain 
		\[-\nabla f(x^*)\in\limsup\limits_{
			x\rightarrow x^*,~
			\alpha\downarrow 0,~\beta\downarrow 0
		}\Theta(x,\alpha,\beta)\subseteq\Theta(x^*,0,0),\]
		so $x^*$ is CC-M-stationary. 
		
		$(\romannumeral3)$~We take $\theta^*\in\limsup\limits_{
			x\rightarrow x^*,~
			\alpha\downarrow 0,~\beta\downarrow 0
		}\Theta(x,\alpha,\beta)$, and define $f=x^T\theta^*$, then $\nabla f(x^*)=\theta^*$. Since
			\[\theta^*\in\limsup\limits_{
				x\rightarrow x^*,~
				\alpha\downarrow 0,~\beta\downarrow 0
			}\Theta(x,\alpha,\beta)\quad\Longrightarrow\quad \theta^*\in \Theta(x^*,0,0),\]
			then $\limsup\limits_{
				x\rightarrow x^*,~
				\alpha\downarrow 0,~\beta\downarrow 0
			}\Theta(x,\alpha,\beta)\subseteq \Theta(x^*,0,0)$.
		\qed
	%\end{proof}	
	
	In Theorem \ref{thm1}, we have proved that any local minimizer of CCOP (or $TNLP(x^*)$) satisfies the sequential optimality condition (CC-PAM-stationarity), and Theorem \ref{PAM->M}$~(\romannumeral2)$ explains	
	\begin{equation}
	\label{SCQ}
	\text{CC-PAM}~+~\text{CC-PAM-regularity}\quad\Longrightarrow\quad \text{CC-M},
	\end{equation}
	in other words, CC-PAM-regularity condition is a CC-CQ. Literature \cite{ALGENCAN} calls the constraint qualification that satisfies the property \eqref{SCQ} as the strict constraint qualification (SCQ). And the conclusion $(\romannumeral3)$ means that the CC-PAM-regularity condition is the weakest SCQ relative to CC-PAM-stationarity. Next, we will apply CC-PAM-stationarity and CC-PAM-regularity condition to enhance the theoretical results of the augmented Lagrangian method.
	
\section{Convergence of Safeguarded Augmented Lagrangian Method}
   	This section will discuss the convergence of using safeguarded augmented Lagrangian method to directly solve the relaxation problem \eqref{relax_problem}. Let $\Lambda=(\lambda,\mu,\gamma,\delta,\eta)\in\cR_+^m\times\cR^p\times\cR^n\times\cR_+\times\cR_+ ^n$, $\rho>0$, then the augmented Lagrangian function of the problem \eqref{relax_problem} can be written as
   		\begin{equation}
   		\begin{aligned}
   		\cL(x,&y,\Lambda,\rho)=f(x)+\frac{\rho}{2}\left[\left\|\left(g(x)+\frac{\lambda}{\rho}\right)_{+}\right\|^{2}+\left\| h(x)+\frac{\mu}{\rho}\right\|^{2}+\right.\\
   		&\left.\left\| x\circ y+\frac{\gamma}{\rho}\right\|^{2}+\left\| \left(n-\kappa-e^{T}y+\frac{\delta}{\rho}\right)_{+}\right\|^{2}+\left\| \left(y-e+\frac{\eta}{\rho}\right)_{+}\right\|^{2}\right].
   		\end{aligned}
   		\label{aug_funtion}
   		\end{equation}
   		Now we give safeguarded augmented Lagrangian method \cite{ALGENCAN}.
   		\begin{alg}
   		\label{SALM}
   		Safeguarded Augmented Lagrangian Method (SALM)\\
   		%\begin{enumerate}[\textbf{Step} \upshape\bfseries1]
   		 {\upshape\bfseries Step 1}(Initialization) Given $(x^{0},y^{0})\in\cR^{n}\times\cR^{n}$, $\lambda_{max}>0$, $\mu_{min}<\mu_{max}$, $\gamma_{min}<\gamma_{max}$, $\delta_{max}>0$, $\eta_{max}>0$, $\tau>1$, $\sigma>1$, $\{\epsilon_k\}\in\cR_+$ and $\epsilon_k\downarrow 0$. Choose initial values $\bar{\lambda}^0\in[0,\lambda_{max}]^m$, $\bar{\mu}^0\in[\mu_{min},\mu_{max}]^p$, $\bar{\gamma}^0\in[\gamma_{min},\gamma_{max}]^n$, $\bar{\delta}^0\in[0,\delta_{max}]$, $\bar{\eta}^0\in[0,\eta_{max}]^n$, $\rho_{0}>0$ and set $k=1$.\\
   		{\upshape\bfseries Step 2}(Update of the iterate) Compute $(x^{k},y^{k})$ as an approximate solution of
   				\begin{equation}
   				\min~\cL(x,y,\bar{\Lambda}^{k-1},\rho_{k-1})
   				\label{AL_sub_pro}
   				\end{equation}
   				such that  
   				\begin{equation}
   				\parallel \nabla\cL(x^k,y^k,\bar{\Lambda}^{k-1},\rho_{k-1})\parallel\leq\epsilon_k.
   				\label{nabla->0}
   				\end{equation}%}			
   		{\upshape\bfseries Step 3}(Update of the approximate multipliers) 
   				\begin{align}
   				\lambda^k&=(\rho_{k-1} g(x^k)+\bar{\lambda}^{k-1})_+,\hspace{3em}\mu^k=\rho_{k-1} h(x^k)+\bar{\mu}^{k-1},\label{lamada_update}\\
   				\gamma^k&=\rho_{k-1} x^k\circ y^k+\bar{\gamma}^{k-1},\hspace{3.5em}\eta^k=(\rho_{k-1} (y^k-e)+\bar{\eta}^{k-1})_+,\\
   				\delta^k&=(\rho_{k-1} (n-\kappa-e^Ty^k)+\bar{\delta}^{k-1})_+.\quad\  \label{eta_update}
   				\end{align}%}
   		{\upshape\bfseries Step 4} (Update of the penalty parameter) Take 
   				\begin{align*} 
   				U^k&:=\min\left\{-g(x^k),\frac{\bar{\lambda}^{k-1}}{\rho_{k-1}}\right\},~V^k:=\min\left\{-(n-\kappa-e^Ty^k),\frac{\bar{\delta}^{k-1}}{\rho_{k-1}}\right\},\\
   				R^k&:=\min\left\{-(y^k-e),\frac{\bar{\eta}^{k-1}}{\rho_{k-1}}\right\},
   				\end{align*}
   				if $k=1$ or 
   				\begin{equation}
   				\begin{aligned}
   				&\max\{\parallel U^{k-1}\parallel,\parallel h(x^{k-1})\parallel,\parallel x^{k-1}\circ y^{k-1}\parallel,\parallel V^{k-1}\parallel,\parallel R^{k-1}\parallel\}\\
   				&\geq \tau\max\{\parallel U^k\parallel,\parallel h(x^k)\parallel,\parallel x^k\circ y^k\parallel,\parallel V^k\parallel,\parallel R^k\parallel\},
   				\end{aligned}
   				\label{rho_bound}
   				\end{equation}
   				set $\rho_k=\rho_{k-1}$, otherwise set $\rho_k=\sigma\rho_{k-1}$.\\
   		{\upshape\bfseries Step 5} (Update of the safeguarded multipliers)
   				Choose $\bar{\lambda}^k\in[0,\lambda_{max}]^m$, $\bar{\mu}^k\in[\mu_{min},\mu_{max}]^p$, $\bar{\gamma}^k\in[\gamma_{min},\gamma_{max}]^n$, $\bar{\delta}^k\in[0,\delta_{max}]$, $\bar{\eta}^k\in[0,\eta_{max}]^n$.

   				Set $k\leftarrow k+1$, go to Step 1.%}
   		%\end{enumerate}
   		\end{alg}

   		Algorithm \ref{SALM} introduces the safeguard multiplier based on the classic augmented Lagrangian method. The convergence is further improved, and the whole sequence convergence required in the classic methods is relaxed to the subsequence convergence \cite{Kanzow_2017}. The updated way of safeguard multiplier in Step 5 is not unique, such as the most popular projection method. In addition, it should be emphasized that the stopping criterion is not set in Algorithm \ref{SALM}, and we will explore this issue in the subsequent convergence analysis.
   		
   		Before proceeding to the analysis of convergence, a useful assumption is given.
   		\begin{assumption}
   			Assuming that $g:\cR^n\rightarrow \cR^m$ and $h:\cR^n\rightarrow \cR^p$ in CCOP are  semialgebraic function.
   			\label{assum1}
   		\end{assumption}
   		
   	In Sect.2, we have introduced the basic concepts and properties of semialgebraic, explaining that Assumption \ref{assum1} is a relatively mild condition, which covers a large class of problems. In the subsequent analysis, we will see that the "semialgebraic" assumption can be further relaxed. Let $p(x,y,\Lambda,\rho)=\frac{1}{\rho}\left[\cL(x,y,\Lambda,\rho)-f(x)\right]$, it is actually the second penalty part of \eqref{aug_funtion} (excluding penalty parameters). Under the assumption $\ref{assum1}$, the following conclusions can be easily obtained by Lemma \ref{semialge}.	
   		\begin{lemma}
   			If Assumption \ref{assum1} holds, then for any given $\Lambda$, $\rho$, $p(x,y,\Lambda,\rho)$ semialgebraic.
   			\label{lemma2}
   		\end{lemma}
   		
   	Let $\{(x^k,y^k)\}$ be the iterative sequence generated by Algorithm \ref{SALM}. We already know that if $\{x^k\}$ is bounded on a subsequence, then $\{y^k\}$ is bounded on the corresponding subsequence (for detailed proof, see \cite{Schwartz_2021_ALA}). This property shows that the boundedness of the entire iteration sequence can be obtained only by ensuring that it is bounded in the $x$ space (that is, in the original problem). A sufficient condition is given below.
   		\begin{lemma}
   			If $f$ is level bounded, then for any given $\Lambda$, $\rho$, $\cL(x,y,\Lambda,\rho)$ is also level bounded.
   			\label{level_bound}
   		\end{lemma}	
   	    {\it Proof}
   	    %\begin{proof}
   	    ~Let
   	    \begin{equation*}
   	    T(y,\delta,\eta,\rho)=\frac{\rho}{2}\left[\left\|\left(n-\kappa-e^{T}y+\frac{\delta}{\rho}\right)_{+}\right\|^{2}+\left\| \left(y-e+\frac{\eta}{\rho}\right)_{+}\right\|^{2}\right].
   	    \end{equation*}
   	    Meanwhile, for any $\alpha\in\cR$, set $f_{\alpha}$, $T_{\alpha}$ as their respective levels under the $\alpha$ level set. By hypothesis, $f_{\alpha}$ is bounded, and for $T(y,\delta,\eta,\rho)$, we have 
   	    \[\parallel y\parallel\rightarrow\infty\quad\Longrightarrow\quad T(y,\delta,\eta,\rho)\rightarrow\infty,\]
   	    therefore, $T_{\alpha}$ is also bounded. Because of
   	    \[\{(x,y):\cL(x,y,\Lambda,\rho)\leq\alpha\}\subseteq f_{\alpha}\times T_{\alpha},\]
   	    then $\cL(x,y,\Lambda,\rho)$ is level bounded.	
   	    \qed

   	    Lemma \ref{level_bound} shows that $\cL$ is  consistent level bounded about $\Lambda$ and $\rho$. If the subproblem \eqref{AL_sub_pro} is solved using a descent algorithm, then the sequence generated by Algorithm \ref{SALM} must be bounded. Let us now discuss the convergence of Algorithm \ref{SALM}.
   	
    \begin{theorem}
	\label{convergence}
	Let $(x^*,y^*)$ be an accumulation point of $\{(x^k,y^k)\}$ generated by Algorithm \ref{SALM}, that is feasible for the problem \eqref{relax_problem}, and Assumption \ref{assum1} holds, then $x^*$ is a CC-PAM-stationary point.	
	\end{theorem}
	{\it Proof}~
	%\begin{proof}
	For simplicity, we assum, $(x^k,y^k)\rightarrow (x^*,y^*)$. By \eqref{nabla->0}, we obtain
	\[\nabla f(x^k)+\nabla g(x^k)\lambda^k+\nabla h(x^k)\mu^k+\gamma^k\circ y^k\rightarrow 0.\]
	Now let's do the proof in two cases.\\	
	$\romannumeral1)~\{\rho_k\}$ is bounded.
	
	If $\{\rho_k\}$ is bounded, combined with $\{\bar{\Lambda}^k\}$ being bounded and \eqref{lamada_update}-\eqref{eta_update}, we know that $\{\Lambda^k\}$ is also bounded. To avoid repeatedly taking subsequence, we assum $\Lambda^k\rightarrow\Lambda$, then
	\[\nabla f(x^*)+\nabla g(x^*)\lambda+\nabla h(x^*)\mu+\gamma\circ y^*=0.\]
	
	Let $\cI=\{i:\lambda_{i}>0\}$, $\cJ=\{j:\mu_{j}\neq 0\}$, $\cK=\{\imath:\gamma_\imath y_{\imath}^*\neq 0\}\subseteq I_0(x^*)$. If the three are all empty sets, then let $\hat{x}^k=x^*$, $\hat{\lambda}^k=\hat{\mu}^k=\hat{\gamma}^k=0$,  it can be concluded that $x^*$ is a CC-PAM-stationary point. Conversely, if there is at least one non-empty, it can be obtained from Lemma 1 of \cite{Andreani_2012}, there exist $I\subseteq\cI$, $J\subseteq\cJ$, $K\subseteq\cK$ and $(\hat{\lambda}_I,\hat{\mu}_J,\hat{\gamma}_K)$ such that
	\begin{align}
	&\nabla f(x^*)+\sum_{i\in I}\hat{\lambda}_i\nabla g_i(x^*)+\sum_{j\in J}\hat{\mu}_j\nabla h_j(x^*)+\sum_{\imath\in K}\hat{\gamma}_\imath e_\imath=0\label{hat-nabla};\\
	&\hat{\lambda}_i\cdot\lambda_i>0,~i\in I;\label{hat-lambda}\\
	&\hat{\mu}_j\cdot\mu_j>0,~j\in J;\label{hat-mu}\\
	&\hat{\gamma}_\imath\cdot(\gamma_\imath y_\imath^*)>0,~\imath\in K;\label{hat-gamma}
	\end{align}
	And the vector group
	\[\cF=\{\nabla g_{i}(x^*)~(i\in I),~\nabla h_{j}(x^*)~(j\in J),~e_{\imath}~(\imath\in K)\}\]
	is linearly independent.
	
	Take
	\begin{equation}
	\hat{\lambda}_i^k=\left\{\begin{aligned}
	&\hat{\lambda}_i\hspace{1.5em}i\in I,\hspace{0.9em}\\
	&0\hspace{0.9em}otherwise,
	\end{aligned}\right.\hspace{0.2em}	\hat{\mu}_j^k=\left\{\begin{aligned}
	&\hat{\mu}_j\hspace{1.5em}j\in J,\hspace{0.9em}\\
	&0\hspace{0.9em}otherwise,
	\end{aligned}\right.\hspace{0.2em}
	\hat{\gamma}_\imath^k=\left\{\begin{aligned}
	&\hat{\gamma}_\imath~\hspace{1.5em}\imath\in K,\hspace{0.75em}\\
	&0\hspace{0.9em}otherwise,
	\end{aligned}\right.
	.\label{hat}
	\end{equation}
	
	The next key problem is to find a sequence $\{\hat{x}^k\}$, $\hat{x}^k\rightarrow x^*$ such that $\{(\hat{x}^k,\hat{\lambda}^k,\hat{\mu}^k,\hat{\gamma}^k)\}$ is a CC-PAM sequence. Let  
	\[J_+:=\{j:\hat{\mu}_j>0\},~J_-:=\{j:\hat{\mu}_j<0\},~K_+:=\{\imath:\hat{\gamma}_\imath>0\},~K_-:\{\imath:\hat{\gamma}_\imath<0\},\]
	obviously $J_+\bigcup J_-=J$, $K_+\bigcup K_-=K$. And we define
	\begin{align*}
	Z:=\{x:g_i(x)\geq 0,&~h_{j_+}(x)\geq 0,~h_{j_-}(x)\leq 0,~x_{\imath_+}\geq 0,~x_{\imath_-}\leq 0,\\
	&i\in I,~j_+\in J_+,~j_-\in J_-,~\imath_+\in K_+,~\imath_-\in K_-\}.
	\end{align*}
	Since the vector group $\cF$ is linearly independent, then $Z$ satisfies LICQ at $x^*$, then MFCQ must also be satisfied. Thus, there exists $d\in\cR^n$ such that
	\begin{equation*}
	\nabla g_i(x^*)^Td>0,~\nabla h_{j_+}(x^*)^Td>0,~\nabla h_{j_-}(x^*)^Td<0,~e_{\imath_+}^Td>0,~e_{\imath_-}^Td< 0,
	\end{equation*}
	where $i\in I,~j_+\in J_+,~j_-\in J_-,~\imath_+\in K_+,~\imath_-\in K_-$. For simplicity, we set $\parallel d\parallel=1$, and take 
	\[\hat{x}^k=x^*+t_{k}d,\]
	where $t_k\downarrow 0$, this implies $\hat{x}^k\rightarrow x^*$. By \eqref{hat-nabla} and \eqref{hat}, we have
	\[\nabla f(\hat{x}^k)+\nabla g(\hat{x}^k)\hat{\lambda}^k+\nabla h(\hat{x}^k)\hat{\mu}^k+\hat{\gamma}^k\rightarrow 0.\]
	
	Let $i\notin I_g(x^*)$, we have $g_i(x^*)<0$, and $g_i(x^k)<0$ when $k$ is sufficiently large. By \eqref{rho_bound}, we know
	\[\parallel U^k\parallel\rightarrow0~\Longrightarrow~\bar{\lambda}_i^{k-1}\rightarrow 0.\]
	Futhermore, \eqref{lamada_update} implies
	\[\lambda_i^k=(\rho_{k-1} g_i(x^k)+\bar{\lambda}_i^{k-1})_+=0,\quad \text{for all k sufficiently large} ,\]
	then $\lambda_i=0$, that is, $i\notin\cI$. Hence, by \eqref{hat}, we have
	\begin{equation}
	\hat{\lambda}_i^k=0,\quad \forall i \notin I_g(x^*).
	\label{b-lam-satis}
	\end{equation}
	
	Take an index $\imath\in I_\pm(x^*)$, we know $x_\imath^*\neq0$, $y_\imath^*=0$, then $\gamma_\imath y_\imath^*=0$, i.e. $\imath\notin\cK$. By \eqref{hat}, we have
	\begin{equation}
	\hat{\gamma}_\imath^k=0, \quad \imath\in I_\pm(x^*).
	\label{b-gamma-satis}
	\end{equation}
	
	The following verifies that $\{(\hat{x}^k,\hat{\lambda}^k,\hat{\mu}^k,\hat{\gamma}^k)\}$ satisfies conditions $(c)$-$(e)$ of Definition \ref{CC-PAM}. Set
	\begin{align*}
	\Gamma:=\min\{\nabla g_i(x^*)^Td,&~\nabla h_{j_+}(x^*)^Td,~-\nabla h_{j_-}(x^*)^Td,~e_{\imath_+}^Td,~-e_{\imath_-}^Td,\\
	&i\in I,~j_+\in J_+,~j_-\in J_-,~\imath_+\in K_+,~\imath_-\in K_-\}.
	\end{align*}
	
	If $\hat{\lambda}_i^k\neq 0$, then $i\in I$. By \eqref{b-lam-satis}, we know $I\subseteq I_g(x^*)$. This implies
	\begin{align*}
	g_i(\hat{x}^k)&=g_i(x^*)+\nabla g_i(x^*)^T(\hat{x}^k-x^*)+r(\parallel \hat{x}^k-x^*\parallel)\\
	&=t_k\nabla g_i(x^*)^Td+r(t_k),
	\end{align*}
	where $r(t_k)$ represents the high-order infinitesimal of $t_k$. Divide both sides of the above formula by $t_k$, when $k$ is sufficiently large, we have
	\[\frac{g_i(\hat{x}^k)}{t_k}=\nabla g_i(x^*)^Td+\frac{r(t_k)}{t_k}\geq\frac{\Gamma}{2}>0,\]
	hence, $\hat{\lambda}_i^kg_i(\hat{x}^k)>0$.
	
	If $\hat{\mu}_j^k\neq 0$, then $j\in J$. We only discuss $j\in J_-$ here (similarly available for $j\in J_+$). For any $j\in J_-$, when $k$ is sufficiently large, we know
	\[\frac{h_j(\hat{x}^k)}{t_k}=\nabla h_j(x^*)^Td+\frac{r(t_k)}{t_k}\leq-\frac{\Gamma}{2}<0,\]
	then $\hat{\mu}_j^kh_j(\hat{x}^k)>0$.
	
	Analogously, if $\hat{\gamma}_\imath^k\neq0$, by \eqref{hat}, we know $\imath\in K$, then 
	\begin{align*}
	x_\imath^k&=x_\imath^*+t_kd_\imath=t_ke_\imath^Td_\imath<0,\quad&\imath\in K_-,\\
	x_\imath^k&=x_\imath^*+t_kd_\imath=t_ke_\imath^Td_\imath>0,\quad&\imath\in K_+,
	\end{align*}
	hence, for any $\hat{\gamma}_\imath^k\neq0$, obviously $\hat{\gamma}_\imath^kx_\imath^k>0$.
	
	In summary, we show that $\{(\hat{x}^k,\hat{\lambda}^k,\hat{\mu}^k,\hat{\gamma}^k)\}$ is a CC-PAM sequence. Therefore $x^*$ is CC-PAM-stationary point.\\
	$\romannumeral2)~\{\rho_k\}$ is unbounded.
	
	Let $x^k\rightarrow x^*$, by \eqref{nabla->0}, we know
	\[\left\|\nabla_{x} \cL\left(x^{k}, y^{k}, \bar{\Lambda}^{k-1}, \rho_{k-1}\right)\right\|=\left\|\nabla f(x^k)+\nabla g\left(x^{k}\right) \lambda^{k}+\nabla h\left(x^{k}\right) \mu^{k}+\gamma^k\circ y^{k}\right\| \leqslant \varepsilon_{k}.\]
	Set $\tilde{\gamma}^k=\gamma^k \circ y^k$, since $\epsilon_k\downarrow 0$, this implies
	\[\nabla f(x^k)+\nabla g\left(x^{k}\right) \lambda^{k}+\nabla h\left(x^{k}\right) \mu^{k}+\tilde{\gamma}^k\rightarrow 0.\]
	
	Let $i\notin I_g(x^*)$, we have $g_i(x^*)<0$, and $g_i(x^k)<0$ when $k$ is sufficiently large. Since $\rho_k\rightarrow\infty$, then
	\[\lambda_{i}^k=(\rho_{k-1} g_i(x^k)+\bar{\lambda}_i^{k-1})_+=0.\]
	
	Take an index $\imath\in I_\pm(x^*)$, we have $x_\imath^*\neq 0$ and $y_\imath^*=0$, then $y_\imath^k\rightarrow 0$. Thus
	\begin{align*}
		\lim _{k \rightarrow \infty} \gamma_{\imath}^{k} y_{\imath}^{k} &=\lim _{k \rightarrow \infty} \rho_{k-1} x_{\imath}^{k}\left(y_{\imath}^{k}\right)^{2}+\lim _{k \rightarrow \infty}\bar{\gamma}_{\imath}^{k-1} y_{\imath}^{k} \\
		&=\lim _{k \rightarrow \infty} \frac{1}{x_{\imath}^{k}} \rho_{k-1}\left(x_{\imath}^{k} y_{\imath}^{k}\right)^{2}.
	\end{align*}
	
	 Now, Let's prove $\rho_{k-1}\left(x_{\imath}^{k} y_{\imath}^{k}\right)^{2}\rightarrow 0$. For simplicity, we abbreviate $p(x,y,\Lambda,\rho)$ in Lemma \ref{lemma2} as $p(x)$, and define 
	\begin{equation*}
	\bar{p}(x)=\frac{1}{2}\left[\left\|\left(g(x)\right)_+\right\|^2+\left\|h(x)\right\|^2+\left\|x\circ y\right\|^2+\left\|\left(n-\kappa-e^Ty\right)_+\right\|^2+\left\|\left(y-e\right)_+\right\|^2\right].
	\end{equation*}
	 Then 
	\begin{align*}
	\rho_{k-1}\nabla_{(x,y)}p(x^k)&-\rho_{k-1}\nabla_{(x,y)}\bar{p}(x^k)\\
	&\leq\left(\begin{array}{c}
	\nabla g\left(x^{k}\right) \bar{\lambda}^{k-1}+\nabla h\left(x^{k}\right) \bar{\mu}^{k-1}+\bar{\gamma}^{k-1} \circ y^{k} \\
	-\bar{\delta}^{k-1} e+\bar{\eta}^{k-1}+\bar{\gamma}^{k-1} \circ x^{k}
	\end{array}\right),
	\end{align*}
	where, the inequality sign comes from the Lipschitz property of $(\cdot)_+$, such as
	\begin{equation*}
	\lambda^k-(\rho_{k-1}g(x^k))_+=(\rho_{k-1}g(x^k)+\bar{\lambda}^{k-1})_+-(\rho_{k-1}g(x^k))_+\leq\bar{\lambda}^{k-1},
	\end{equation*}
	the others are similar. Since $\bar{\Lambda}^{k-1}$ is bounded, $g,~h\in C^1$ and $(x^k,y^k)\rightarrow(x^*,y^*)$, then there exists $M_1>0$ such that
	\[\rho_{k-1}\nabla_{(x,y)}p(x^k)-\rho_{k-1}\nabla_{(x,y)}\bar{p}(x^k)\leq M_1.\]
	
	On the other hand, by \eqref{nabla->0}, we know
	\begin{align*}
	\left\| \rho_{k-1}\nabla_{(x,y)}p(x^k)\right\|&\leq \left\|\nabla \cL\left(x^{k}, y^{k}, \bar{\Lambda}^{k-1}, \rho_{k-1}\right)\right\|+\left\| \nabla f(x^k)\right\|\\
	&\leq \epsilon_k+\left\|\nabla f(x^k)\right\|.
	\end{align*}
	Thus
	\begin{align*}
	\left\|\rho_{k-1}\nabla_{(x,y)}\bar{p}(x^k)\right\|&=\left\|\rho_{k-1}\nabla_{(x,y)}p(x^k)-(\rho_{k-1}\nabla_{(x,y)}p(x^k)-\rho_{k-1}\nabla_{(x,y)}\bar{p}(x^k))\right\|\\
	&\leq\epsilon_k+\left\|\nabla f(x^k)\right\|+M_1.
	\end{align*}
	Furthermore, since $\epsilon_k\downarrow 0$, $f\in C^1$ and $x^k\rightarrow x^*$, then there exists $M>0$ such that
	\begin{equation}
	\left\|\rho_{k-1}\nabla_{(x,y)}\bar{p}(x^k)\right\|\leq M\label{semi1}.
	\end{equation}
	
	At the same time, by Assumption \ref{assum1} and Lemma \ref{semialge}, we know $\bar{p}(x)$ is semialgebraic. So there exists $C>0$, $\theta\in[0,1)$ such that
	\begin{equation}
	\left\|\rho_{k-1}\nabla_{(x,y)}\bar{p}(x^k)\right\|\geq C\rho_{k-1}\left\|\bar{p}(x^k)-\bar{p}(x^*)\right\|^\theta=C\rho_{k-1}\left\|\bar{p}(x^k)\right\|^\theta.
	\label{semi2}
	\end{equation}
	By \eqref{semi1} and \eqref{semi2}, we obtain
    \begin{align*}
    0\leq\rho_{k-1}\left\|x^k\circ y^k\right\|^2\leq\left\|\rho_{k-1}\bar{p}(x^k)\right\|&\leq C^{-1}\left\|\rho_{k-1}\nabla_{(x,y)}\bar{p}(x^k)\right\|\left\|\bar{p}(x^k)\right\|^{1-\theta}\\
    &\leq C^{-1}M\left\|\bar{p}(x^k)\right\|^{1-\theta}\rightarrow 0.
    \end{align*}
    Thus
    \[\lim\limits_{k\rightarrow\infty}\rho_{k-1}x_\imath^k(y_\imath^k)^2=\lim\limits_{k\rightarrow\infty}\frac{1}{x_\imath^k}\rho_{k-1}(x_\imath^ky_\imath^k)^2=0,\]
    that is $\lim\limits_{k\rightarrow\infty}\tilde{\gamma}_\imath^k=0$. Set
    \begin{equation}
    \hat{\gamma}_\imath^k=\left\{\begin{aligned}
    &0&\imath\in I_\pm(x^*)\\
    &\tilde{\gamma}_\imath^k&\imath\in I_0(x^*)
    \end{aligned}\right.,
    \label{hat-gamma2}
    \end{equation}
	It's easy to know that there are still
	\[\nabla f(x^k)+\nabla g\left(x^{k}\right) \lambda^{k}+\nabla h\left(x^{k}\right) \mu^{k}+\hat{\gamma}^k\rightarrow 0.\]
	
	Let $\pi_k:=\left\|\left(1,\lambda^k,\mu^k,\hat{\gamma}^k\right)\right\|_\infty$. If $\{\pi_k\}$ is bounded, the proof process of case $\romannumeral1)$ can verify that the conditions $(c)$-$(e)$ are established. Now, we consider the case that $\pi_k$ is unbounded. If $\lim\limits_{k \rightarrow \infty}\frac{\lambda_{i}^k}{\pi_k}>0$, then
	\[\lim_{k \rightarrow \infty}\frac{\lambda_{i}^k}{\pi_k}=\lim_{k \rightarrow \infty}\frac{\rho_{k-1} g_i(x^k)+\bar{\lambda}_i^{k-1}}{\pi_k}=\lim_{k \rightarrow \infty}\frac{\rho_{k-1} g_i(x^k)}{\pi_k}>0.\]
	Obviously, we have $g_i(x^k)>0$ for all $k$ sufficiently
	large, so
	\[\lambda_{i}^k g_i(x^k)>0,\quad if ~\lim\limits_{k \rightarrow \infty}\frac{\lambda_{i}^k}{\pi_k}>0.\]
	
	If $\lim\limits_{k \rightarrow \infty}\frac{|\mu_{j}^k|}{\pi_k}>0$, then
	\[\lim_{k \rightarrow \infty}\frac{\mu_{j}^k}{\pi_k}=\lim_{k \rightarrow \infty}\frac{\rho_{k-1} h_j(x^k)+\bar{\mu}_j^{k-1}}{\pi_k}=\lim_{k \rightarrow \infty}\frac{\rho_{k-1} h_j(x^k)}{\pi_k}.\]
	Observe that $\mu_{j}^k$ has the same sign as $h_j(x^k)$, this implies
	\[\mu_{j}^k h_j(x^k)>0,\quad if~\lim\limits_{k \rightarrow \infty}\frac{|\mu_{j}^k|}{\pi_k}>0.\]
	
	Similarly, if $\lim\limits_{k \rightarrow \infty}\frac{|\hat{\gamma}_{\imath}^k|}{\pi_k}>0$, by \eqref{hat-gamma2}, we know $\imath\in I_0(x^*)$, then $y_\imath^k\rightarrow y_\imath^*\neq 0$. Meanwhile
	\begin{equation*}
	\lim_{k \rightarrow \infty}\frac{\hat{\gamma}_{\imath}^k}{\pi_k}=\lim_{k \rightarrow \infty}\frac{\gamma_\imath^k y_\imath^k}{\pi_k}=\lim_{k \rightarrow \infty}\frac{\rho_{k-1} x_\imath^k(y_\imath^k)^2}{\pi_k}.
	\end{equation*}
	Therefore, when $k$ is sufficiently large, $\hat{\gamma}_\imath^k$ has the same sign as $x_\imath^k$, namely
	\[\hat{\gamma}_\imath^k x_\imath^k>0,\quad if ~\lim\limits_{k \rightarrow \infty}\frac{|\hat{\gamma}_{\imath}^k|}{\pi_k}>0.\]
We have completed this proof.\qed
    %\end{proof}
    
    As can be seen from the proof of  Theorem \ref{convergence}, the "semi-algebraic" hypothesis is essentially to ensure that for any given $\Lambda$, $\rho$, $p(x,y,\Lambda,\rho)$ has KL properties. Therefore, Assumption \ref{assum1} can be further relaxed to the structure of $\cO$-$minimal$, or even to the assumption that $p(x,y,\Lambda,\rho)$ has KL properties (that is, the same as the conditions of \cite{Schwartz_2021_Sequ}), the conclusion of Theorem \ref{convergence} still holds. There are two reasons why we did not do this. One is that $y$ is an artificial variable, so all assumptions should not be imposed on the $y$ space; in addition, through Lemma \ref{lemma2} and Lemma \ref {level_bound} can show that the introduction of $y$ does not destroy the good properties of the hypothesis on the $x$ space.
	
	Theorem \ref{convergence} states that any feasible accumulation point of Algorithm \ref{SALM} is a CC-PAM-stationary point if Assumption \ref{assum1} holds. And from the proof process, it can be seen that the sequence generated by Algorithm \ref{SALM} is not necessarily a CC-PAM sequence.  In fact, this is not contradictory, because we require the existence of the corresponding CC-PAM sequence in Definition \ref{CC-PAM}. But at least this shows that it is not appropriate to take CC-PAM-stationarity as the stop criterion. On the other hand, according to Theorem \ref{PAM->M}, if there is an additional CC-PAM regularity condition holds at $x^*$, then it is a CC-M-stationary point. In other words, when CC-PAM regularity condition is established, CC-M-stationarity itself is a very suitable stop criterion, \cite{Schwartz_2021_ALA} has verified the validity of this method, this paper mainly emphasizes the theoretical improvement, not to repeat the experiment. In addition, in conjunction with Proposition \ref{prop2}, it is clear that the following conclusion holds.
	\begin{theorem}
			Let $(x^*,y^*)$ be an accumulation point of $\{(x^k,y^k)\}$ generated by Algorithm \ref{SALM}, Assumption \ref{assum1} holds, that $(x^*,y^*)$ is feasible for the relaxation problem \eqref{relax_problem}, and meet CC-PAM regularity condition at $x^*$. Then $(x^*,y^*)$ is a CC-M-stationary point and there exists $z^*\in\cR^n$ such that $(x^*,z^*)$ is a CC-S-stationary point.
	\end{theorem}
   	
 \section{Final Remarks}
    In this paper, we study the continuous relaxation form of CCOP, which is more popular in recent years. We propose a new sequential optimality condition called CC-PAM-stationarity. In Sect.3, we prove that CC-PAM-stationarity is strictly superior to CC-AM-stationarity. Moreover, any local minimizer of CCOP is a CC-PAM-stationary point without any additional assumptions. Obviously, CC-PAM-stationarity is a better measure of optimality than CC-AM-stationarity. In addition, we introduced a new constraint qualification called CC-PAM-regularity in Sect.4, which is weaker than CC-AM-regularity. It is proved that if the CC-PAM regularity condition is established, then any CC-PAM-stationary point all are CC-M-stationary points.
   	
   	In Sect. 5, we apply the new sequential optimality condition proposed in this paper, CC-PAM-stationarity, to the safeguarded augmented Lagrangian method (Algorithm \ref{SALM}), which further improves the existing theoretical results. We have proved that under mild conditions such as KL properties, any feasible convergence point of Algorithm \ref{SALM} is a CC-PAM-stationary point; further, if the CC-PAM-regularity condition, it can converge to a CC-M-stationary (essentially CC-S-stationary) point. In other words, in this case, the CC-M-stationarity is the natural termination criterion of Algorithm \ref{SALM}. Meanwhile, we emphasize that if the same assumptions as the existing results are used, the conclusions of this article are still valid.

\begin{acknowledgements}
The research were partially supported by the Natural Science Foundation of
Shandong Province with No.ZR2019BA014, The Key Research and Development Projects of
Shandong Province with No.2019GGX104089 and the Foundation of Liaoning Educational Committee (LQ2019019).
\end{acknowledgements}

%References
% BibTeX users  please use  \bibliographystyle{spmpsci_unsrt}. The option spmpsci_unsrt prints the references in JOTA format  in the order they are cited.  
%Otherwise, please use the following:
 
%\bibliographystyle{spmpsci_unsrt}	
\bibliographystyle{spmpsci}
\bibliography{refence}
%\begin{thebibliography}{}
%
%% References must be listed in the order in which they actually appear in the text, not alphabetically. 
%
%% Format for journal articles
%\bibitem{1.} Huang, H.Y.: Unified approach to quadratically convergent algorithms for function minimization. J. Optim. Theory Appl. 5, 354-405 (2009)
%
%% Format for reports
%\bibitem{2.} Yang, T. L.: Optimal control for a rocket in a three-dimensional central force field. Technical Memorandum TM 69-2011-2, Bellcomm (1969)      
%                       
%%Format for books
%\bibitem{3.}Nocedal, J., Wright, S.J.: Numerical Optimization, Springer, New York  (2000)
%
%%Format for edited books
%\bibitem{4.} Miele, A. (ed.): Theory of Optimum Aerodynamic Shapes. Academic Press, New York (1965)
%
%%Format for articles in edited books
%\bibitem{5.} Ralston, A.: Numerical integration methods for the solution of ordinary differential equations. In: Ralston, A., Wilf, H.S. (eds.): Mathematical Methods for Digital Computers, vol. 1, pp. 95Ð109. Wiley, New York (1960)
%
%\end{thebibliography}

\end{document}